\documentclass[twoside,leqno,twocolumn]{article}

\usepackage{xurl}

\usepackage{amsmath,amssymb,graphicx,ltexpprt,hyperref,cleveref}
\usepackage[table]{xcolor}
\usepackage{algorithm}

\newtheorem{example}{Example}[section]

\begin{document}

\title{\Large A Uniform Sampling Procedure for Abstract Triangulations of Surfaces}
\author{Rajan Shankar\thanks{School of Mathematics and Statistics F07, The University of Sydney, NSW 2006 Australia, \texttt{rsha5949@uni.sydney.edu.au}}
\and Jonathan Spreer\thanks{School of Mathematics and Statistics F07, The University of Sydney, NSW 2006 Australia, \texttt{jonathan.spreer@sydney.edu.au}}}

\date{}

\maketitle

\fancyfoot[R]{\scriptsize{Copyright \textcopyright\ 2023 by SIAM\\
Unauthorized reproduction of this article is prohibited}}

\begin{abstract} \small We present a procedure to sample uniformly from the set of combinatorial isomorphism types of balanced triangulations of surfaces --- also known as graph-encoded surfaces. For a given number $n$, the sample is a weighted set of graph-encoded surfaces with $2n$ triangles.

  The sampling procedure relies on connections between graph-encoded surfaces and permutations, and basic properties of the symmetric group.
  
  We implement our method and present a number of experimental findings based on the analysis of $138$ million runs of our sampling procedure, producing graph-encoded surfaces with up to $280$ triangles. 
  
  Namely, we determine that, for $n$ fixed, the empirical mean genus $\bar{g}(n)$ of our sample is very close to $\bar{g}(n) = \frac{n-1}{2} - (16.98n -110.61)^{1/4}$. Moreover, we present experimental evidence that the associated genus distribution more and more concentrates on a vanishing portion of all possible genera as $n$ tends to infinity. Finally, we observe from our data that the mean number of non-trivial symmetries of a uniformly chosen graph encoding of a surface decays to zero at a rate super-exponential in $n$.
\end{abstract}

\section{Introduction and Motivation}
 
    The study of surfaces, or, more generally, manifolds, is central to fields of research such as computer graphics, robot motion planning, or data science, and many areas within pure and applied mathematics. Often, these surfaces or manifolds are represented as {\em triangulations}, i.e., as a decomposition into triangles and their higher dimensional analogues, called {\em simplices}. See \Cref{surfaces,triangulations} for examples of surfaces and their triangulations.
    
    Triangulations are often used in an abstract setting, that is, without considering them to be embedded in an ambient space. In this abstract setting it is natural to consider them up to their combinatorial isomorphism types, i.e., up to a relabelling of their simplices. 
    
    It is also common to consider manifolds or surfaces up to continuous deformations, that is, up to {\em homeomorphy}. The famous classification theorem of orientable surfaces then states that for every non-negative integer $g$, there is exactly one such {\em topological type} of surface --- and $g$ is called the {\em genus} of this type. 
        
    Naturally, given positive integers $d$ and $m$, there are many triangulations of $d$-dimensional manifolds with $m$ top-dimensional simplices. But very little is known about this set of triangulations as a whole. 
    
    In this article, we address this issue in the case $d=2$: we present a method that, for $n$ given, uniformly samples balanced triangulations (aka. graph encodings) of surfaces from the list of combinatorial isomorphism types of all such triangulated surfaces with $2n$ triangles. Our sampling procedure is practical for triangulations up to $\sim 200$ triangles and provides novel insights into the set of all graph encodings of surfaces.
  
  There exists a substantial body of literature on models to generate random triangulations of surfaces. Gamburd and Makover \cite{Gamburd-makover} use cubic graphs on Riemannian surfaces to obtain what they call a {\em Random Riemann Surface}. Pippenger and Schleich in \cite{pippenger-schleich} randomly glue triangles along their edges. 
  More recently, Guth, Parlier and Young \cite{guth-parlier-young} use the moduli space of hyperbolic metrics to investigate expected properties of random hyperbolic surfaces. We point to \cite[Section 2]{pippenger-schleich} for a more comprehensive overview of further random surface models.
  
  The work presented here exhibits two key differences: {\em (a)} We entirely focus on balanced triangulations. While we would like to eventually remove this caveat, we still obtain a complete method for one of the standard classes of triangulations presented in literature. {\em (b)} Unlike in the other approaches described above, our method samples surface triangulations uniformly from the list of combinatorial isomorphism types of all triangulations. Combinatorial isomorphy is what computational methods on triangulations see. Uniform sampling up to combinatorial isomorphy is thus the most precise method to judge how algorithms on (abstract) triangulations perform on average, or - equivalently - how complicated triangulations can be with respect to applications. See the following paragraph for examples.

  Despite these key differences, our experimental results are mostly in-line with existing results on the genus distribution, and connectedness of triangulations (cf., for instance, \cite{pippenger-schleich}). Hence, asymptotically, our model behaves very similarly to the other models in the literature. The question of whether we can make this statement more precise is an interesting one which we will not investigate further in this work.

  \paragraph*{Higher dimensions.} The challenge remains to expand this technique to sampling triangulations of $d$-dimensional manifolds, $d\geq 3$, or triangulations of a fixed manifold $M$\footnote{For the latter problem, methods from statistical physics applied to local combinatorial modifications of triangulations (e.g., {\em Pachner moves} \cite{Pachner87KonstrMethKombHomeo}) may offer a viable approach, but serious complications  remain. This is the subject of ongoing research.}
    
    Having knowledge about the properties of random triangulations of manifolds is extremely valuable. In $3$-manifold topology, for instance, many fundamental topological problems have algorithmic solutions \cite{Haken61Unknot, Jaco-algorithm-1984,Rubinstein953SphereRec}. Some of these algorithms have worst-case exponential running times but exhibit polynomial time behaviour \cite{Burton12Unknot} (much like in the case of the celebrated simplex method, \cite{Spielman2009}). A uniform sampling procedure could explain this kind of behaviour and give data on expected running times of important algorithms in the field. 
    
    One fundamental difficulty in defining such a uniform sampling procedure in higher dimensions is the following: Given a set of $m$ simplices of dimension $d \geq 3$, randomly gluing them along their $(d-1)$-dimensional faces (the generalisation of the approach used in \cite{pippenger-schleich}) produces a triangulation $T$ of some $d$-dimensional space. The probability that this space is a manifold decreases to zero, as $m$ tends to infinity. Moreover, checking if a particular choice of $T$ indeed triangulates a manifold is non-trivial for $d=3$, requires worst-case exponential running time in dimension $4$ \cite{Rubinstein953SphereRec}, and is known to be impossible in general for dimensions $>5$ \cite{Markov58HomeoProb,Stillwell93Undecidability}. 

    \paragraph*{Our contribution.} Here, we present a sampling method for dimension two, where the complication described above conveniently disappears: given an even number of triangles, randomly gluing them along their edges always results in a (possibly disconnected) triangulation of a surface.

	To address the still remaining issue of relative probabilities of samples in the space of combinatorial isomorphism types of surfaces, we restrict ourselves to orientable {\em graph-encoded manifolds}, first described in \cite{Ferri76,Pezzana74Crystallizations} (see also the survey \cite{Ferri86Survey}). These can be described by a set of $d$ permutations, where $d$ denotes the dimension. Graph-encoded manifolds can be represented by a graph whose nodes represent simplices and whose coloured arcs describe how simplices are glued together. They are regular simplicial cell-complexes, meaning that their simplices do not admit any self-identifications (in contrast to {\em semi-simplicial triangulations}, where this is possible), but a pair of simplices can intersect in more than one common face (in contrast to {\em simplicial complexes}).

\medskip

In very rough terms, our sampling method {\em (a)} defines a standard way of presenting graph-encoded surfaces; {\em (b)} samples a pair of permutations into this presentation; and {\em (c)} uses algebraic properties of the group of permutations to compute the number of isomorphic copies of the current sample in the space of all possible samples. The inverse of this number is the weight with which the isomorphism type of the graph encoding is added to the sample. 
	
A Python implementation \cite{PythonCoreTeam2020} of our method exhibits near-linear running time in $n$ for obtaining one $2n$-triangle sample. We run the procedure $1$ million times for each value $3 \leq n \leq 140$.
	
	The set of samples experimentally answers questions about the space of surface triangulations. In this work, we focus on the two following questions:	
	\begin{enumerate}
	  \item {\em How is the genus distributed in a uniform sample of graph-encoded surfaces?} For graph encodings with $2n$ triangles, the genus of a sample must lie between $0$ and $\lfloor \frac{n-1}{2}\rfloor $. In our experiments the empirical mean genus is roughly $\bar{g} (n) = \frac{n-1}{2} - (16.98n -110.61)^{1/4}$ and the bulk of our samples has a genus very close to this value. Here, the exponent of the correction term is not claimed to be accurate. Instead, we point to clear experimental evidence that it is strictly between $0$ and $1$.
	  \item {\em What is the average number of symmetries of graph-encoded surfaces?}  The average number of non-trivial symmetries seems to decay to zero at a rate super-exponentially in the number of triangles. Note that a graph-encoded surface with $2n$ triangles can have up to $12n$ symmetries.
	\end{enumerate}

\section{Preliminaries} \label{background}

\paragraph*{Surfaces.} A \textit{surface} is a manifold of dimension two. More precisely, it is a space in which every point has a neighbourhood that is {\em homeomorphic} to $\mathbb{R}^2$. Here, homeomorphic means ``up to continuous deformation''. 

We restrict our research to {\em topological} surfaces that are both \textit{closed} and \textit{orientable}. {\em Topological} means that we consider two surfaces to be the same if one can be turned into the other by a continuous motion (i.e., if they are homeomorphic). This has been proven to be a highly-effective modelling choice in many settings. A {\em closed} surface does not have a boundary, i.e., we cannot ``fall off'' the surface. In an {\em orientable} surface it is not possible to end up on the other side of the surface by simply walking along it. The M\"obius strip is the classic example of a non-orientable surface (with boundary).

\begin{figure}[htb]
	\centerline{\includegraphics[width=0.48\textwidth]{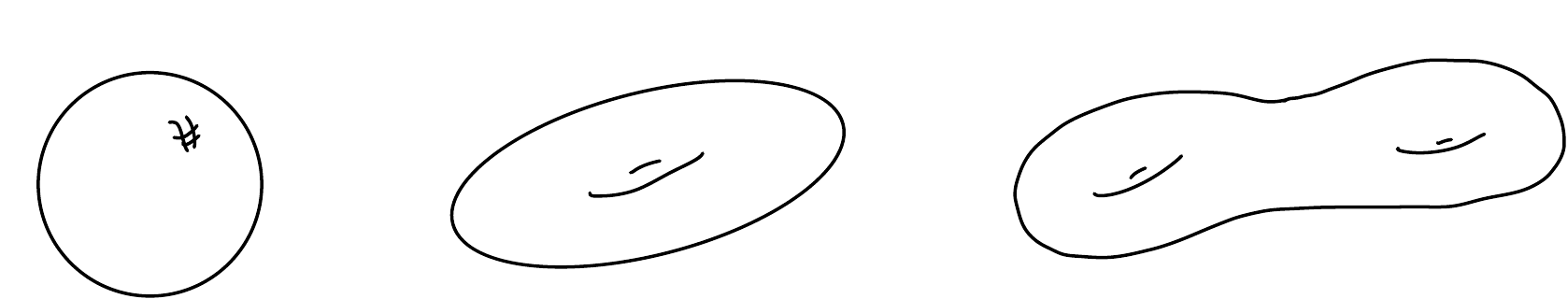}}
	\vspace{-.3cm}
	\caption{The topological closed orientable surfaces of genus zero, one, and two: the $2$-dimensional sphere, the $2$-dimensional torus, and the Pretzel surface.}
	\label{surfaces}
\end{figure}

The {\em genus} of a closed orientable surface is an integer denoting the number of ``holes'' in it. Equivalently, it is the maximum number of cuttings along non-intersecting loops on the surface such that the resulting surface with boundary is still connected. The famous {\em classification of surfaces} states that there is exactly one topological closed orientable surface for every genus $g \geq 0$. The sphere is the surface of genus $0$, the torus is the surface of genus $1$, and so on. See \Cref{surfaces} for examples.

\begin{figure}[htb]
	\centerline{\includegraphics[width=0.48\textwidth]{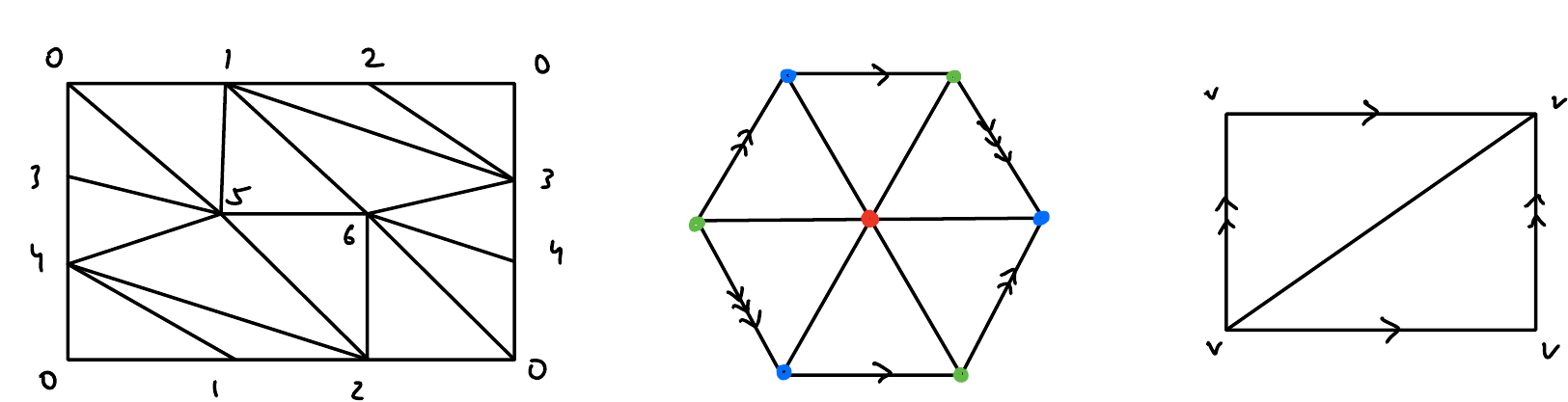}}
	\vspace{-.3cm}
	\caption{Three types of triangulations of the $2$-dimensional torus. For each type, the respective triangulation has the smallest possible number of triangles. Left: simplicial complex with $14$ triangles. Centre: graph-encoded surface with $6$ triangles. Right: semi-simplicial triangulation with $2$ triangles.}
	\label{triangulations}
\end{figure}

\paragraph*{Triangulations and graph encodings.} An \textit{(abstract) triangulation} of a closed surface consists of $2n$ triangles glued together pairwise along their $3\times 2n=6n$ edges such that the resulting space is connected. Note that the pairwise gluings of the edges forces the number of triangles of a triangulated surface to be even. \textit{Abstract} means that the triangulation is not embedded in a space, i.e. the vertices of the triangulation do not have coordinates. See \Cref{triangulations} for three triangulations of the $2$-dimensional torus.

A {\em balanced triangulation of a surface} is a triangulation admitting a {\em rainbow colouring} of its vertices. This means that in a triangulation $T$ with $2n$ triangles we can colour the vertices in three colours such that the vertices of each triangle are coloured in three distinct colours. Throughout this article, we refer to these three colours as the first, the second, and the third colour and draw them in red, green, and blue respectively. 

Given such a colouring of the vertices, we can encode $T$ as a graph $G = (V(G),E(G))$ with an arc colouring $C: E(G) \to \{\text{ red},\text{ green},\text{blue }\}$ through the following procedure: For each triangle in $T$, add one node to $V(G)$. For two triangles of $T$ sharing an edge $e$, add an arc to $E(G)$. Colour this arc in the colour of the unique vertex of the triangles opposite the edge $e$. It follows that $|V(G)|=2n$, and $|E(G)| = 3n$. We call (a drawing of) this arc-coloured graph a \textit{graph-encoded surface} or \emph{graph encoding of a surface}. Note that balanced triangulations and graph encodings are different presentations of the identical combinatorial object: an abstract triangulation of a surface. In this article we mostly work with graph encodings since this representation is more convenient for our purposes.

\begin{figure}[htb]
    \begin{center}
        \centerline{\includegraphics[width=0.48\textwidth]{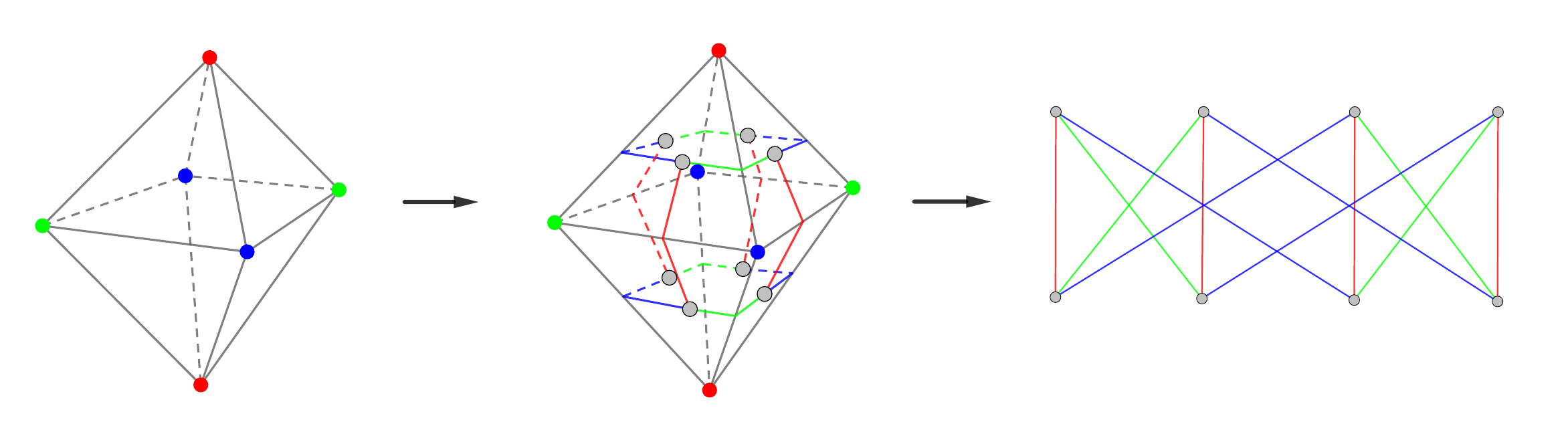}}
    \end{center}
    \vspace{-1.3cm}
    \caption{Graph encoding of the octahedron. Left: vertices are rainbow coloured. Centre: coloured arcs connect centre-points of triangles. Right: rearranging the arcs give a simpler picture.\label{gem_from_triangulation}}
\end{figure}

Given a graph encoding $(G,C)$, $G$ must be bipartite if and only if $(G,C)$ represents an orientable surface \cite{Cavicchioli80Bipartite}. Moreover, the two vertex sets defining the bipartiteness must both be of cardinality $n$. It follows that $(G,C)$ can be re-drawn with arcs of the first colour in vertical lines next to each other, and the remaining two colours connect pairs of nodes from top to bottom. See \Cref{gem_from_triangulation} for an example where $T$ is the octahedron. 

Since $(G,C)$ is an arc-coloured graph with every node adjacent to exactly three arcs, removing all arcs of any one colour leaves a collection of bi-coloured cycles. For instance, in the case of the octahedron, removing arcs of one colour results in two bi-coloured cycles, each of length four. Moreover, these bi-coloured cycles are in 1-to-1 correspondence with the vertices of the underlying triangulation: every bi-coloured cycle of the first and second colour runs around a vertex of the third colour, and so on; see \Cref{gem_from_triangulation}.

We call two graph-encoded surfaces $(G_1,C_1)$ and $(G_2,C_2)$ {\em isomorphic}, if there exist bijections $\phi : V(G_1) \to V(G_2)$ and $\psi: C_1(E(G_1)) \to C_2(E(G_2))$ such that $\phi$ maps arcs of $G_1$ to arcs of $G_2$, and $\psi(C_1 (u,v)) = C_2 (\phi(u),\phi(v))$ for all arcs $(u,v) \in E(G_1)$. If $C_1 = C_2$ and $\psi$ is the identity, $(\phi,\psi)$ is called a {\em colour-preserving isomorphism}. If $G_1 = G_2$ and $\phi$ is the identity, $(\phi,\psi)$ is called a {\em colour swap isomorphism}.

Closely related to isomorphisms are {\em symmetries} of graph encodings. A symmetry of a graph encoding $(G,C)$ is a pair of bijections $\phi : V(G) \to V(G)$ and $\psi: C(E(G)) \to C(E(G))$ such that $(\phi(u),\phi(v)) \in E(G)$ and $C (\phi(u),\phi(v)) = C((u,v))$ for all arcs $(u,v) \in E(G)$. In other words, a symmetry is a reshuffling of the nodes and colours of the graph encoding leaving the graph encoding as a whole unchanged. A symmetry is called {\em colour-preserving} if $\psi$ is the identity, otherwise it is called a {\em colour swap}. Note that an isomorphism can be neither colour-preserving nor a colour swap, and that the only colour-preserving colour swap isomorphism is the identity.

As already exhibited, we reserve the terms ``vertex'' and ``edge'' to be used when discussing triangulations, and the terms ``node'' and ``arc'' when discussing graphs.

\medskip 

\begin{figure}[htb]
    \begin{center}
        \centerline{\includegraphics[width=0.48\textwidth]{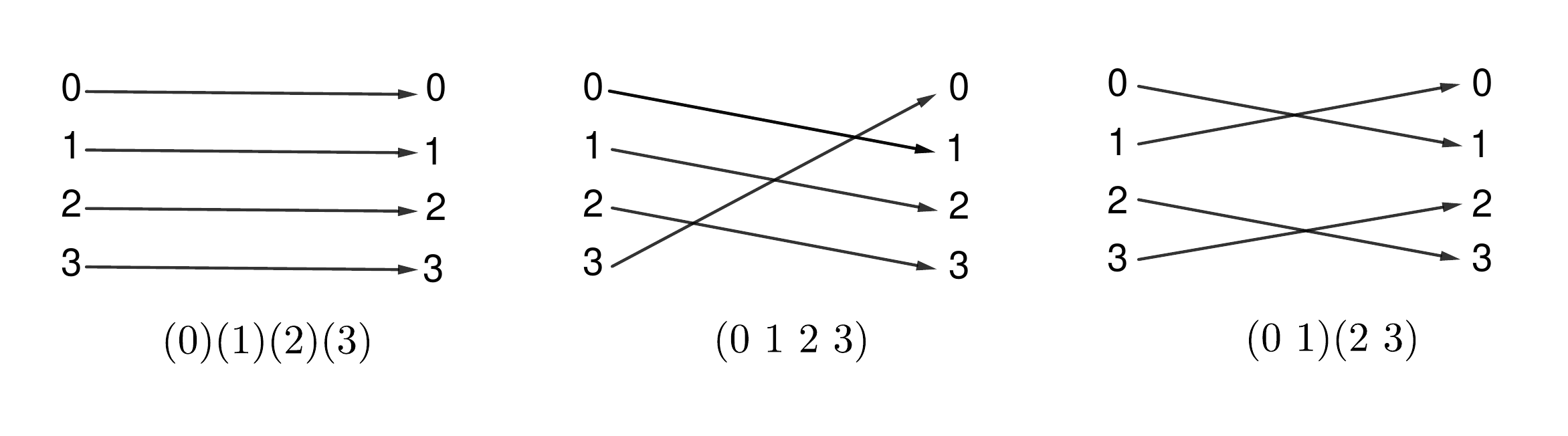}}
    \end{center}
    \vspace{-1.3cm}
    \caption{Three permutations and their cycle notations.}
    \label{cycle_notation}
\end{figure}

\paragraph*{Permutations.} Let $[n] = \{0,1,\dots, n-1\}$ be the set of the first $n$ integers. A \textit{permutation} is a bijection $\sigma : [n] \to [n]$, where $n>0$. Permutations can be written as a vector $(\sigma(0), \sigma(1), \ldots , \sigma(n-1))$ or in {\em cycle notation} $(0, \sigma(0), \sigma^2 (0), \ldots , \sigma^{k-1} (0))(m, \sigma(m), \ldots)\ldots$ where $k$ is the smallest positive integer such that $\sigma^k (0) = 0$ and $m$ is the smallest positive integer not in $\{0, \sigma(0), \sigma^2 (0), \ldots , \sigma^{k-1} (0) \}$. See \Cref{cycle_notation} for examples. 

\begin{figure}[htb]
    \begin{center}
        \centerline{\includegraphics[width=.35\textwidth]{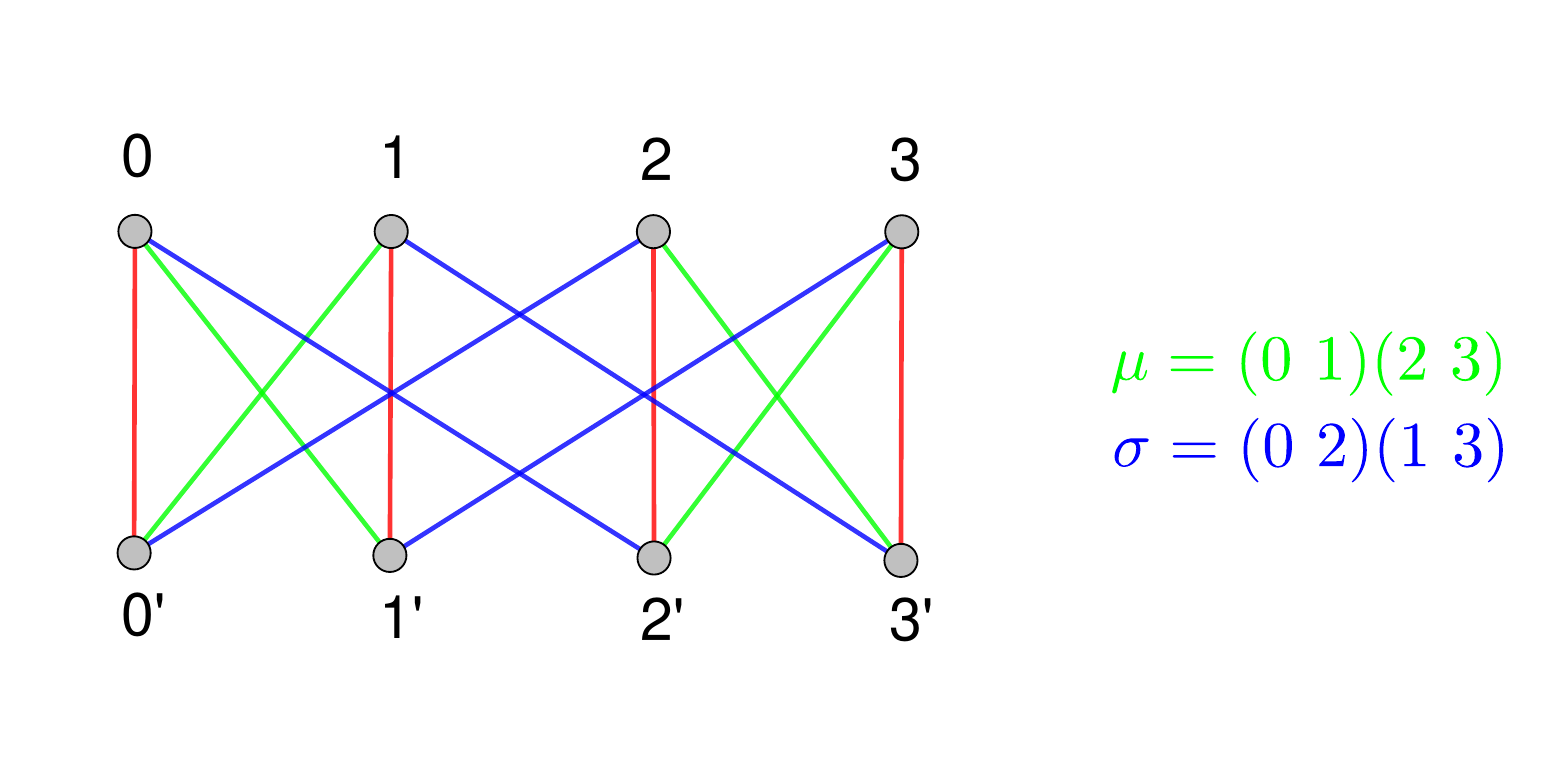}}
    \end{center}
    \vspace{-1.3cm}
    \caption{Converting a graph-encoded surface into a pair of permutations.}
    \label{gem_permutations}
\end{figure}

Looking at the rightmost picture of \Cref{gem_from_triangulation} and numbering the vertices in the top row by $0,1,\ldots ,n-1$ and in the bottom row by $0',1',\ldots ,(n-1)'$, we can interpret each colour set of arcs as a permutation of the nodes in the top row to the nodes in the bottom row. Since we always draw the graph-encoded surface such that one colour appears as vertical arcs (or, the identity permutation), a graph-encoded surface can be encoded by a pair of permutations $(\mu, \sigma)$ determining how the second and third colours are drawn. See \Cref{gem_permutations} for the two permutations encoding the octahedron.

Let $\sigma$ be a permutation of $n$ elements given in cycle notation. Let $p_i$, $1 \leq i \leq k$, be the different lengths of the cycles of $\sigma$, and let $a_i>0$ be the number of cycles of length $p_i$. The set $\{ p_i : a_i \mid 1 \leq i \leq k \}$ is called the {\em cycle structure} of $\sigma$. Naturally, we have $\sum_i a_i \cdot p_i = n$. Moreover, writing down the $p_i$ with multiplicity defines a {\em partition} of $n$, that is, a set of positive integers summing up to $n$. Our convention here is that we write the elements of a partition in weakly decreasing order. See the following section on the symmetric group for some examples of partitions.

Finally, take a minute to verify that every cycle of $\mu$ corresponds to a bi-coloured cycle of the first and second colour, and every cycle of $\sigma$ corresponds to a bi-coloured cycle of the first and third colour. Also note that the bi-coloured cycles of the second and third colour are in 1-to-1 correspondence with the cycles of $\mu^{-1}\circ \sigma$, where $\mu^{-1}$ denotes the inverse of $\mu$ as a bijective function. In particular, there is a 1-to-1 correspondence between the cycles of $\mu$, $\sigma$, and $\mu^{-1}\circ \sigma$ and the vertices of the underlying triangulation.

\medskip

\paragraph*{The symmetric group.} The set of permutations on $n$ elements together with their composition ``$\circ$'' form the {\em symmetric group on $n$ elements} $S_n$. There are $n!$ permutations on $n$ elements and we have for the order of the symmetric group $|S_n| = n!$. 

Given $\sigma \in S_n$, the subset $C_{\sigma} = \{ \mu \circ \sigma \circ \mu^{-1} \mid\, \mu \in S_n \} \subset S_n$ is called the conjugacy class of $\sigma$. Here, $\sigma$ is called the {\em representative} of $C_{\sigma}$. The conjugacy classes of $S_n$ are independent of their representatives, and partition the set of all permutations $S_n$. Moreover, $C_{\sigma}$ contains exactly all the permutations of $\sigma$ with the same cycle structure as $\sigma$. It follows that the conjugacy classes of $S_n$ are in $1$-to-$1$-correspondence with the integer partitions of $n$, where each partition $\lambda$ is associated to the conjugacy class containing the permutations with cycle structure corresponding to $\lambda$. It is hence well-defined to label the conjugacy classes of $S_n$ by partitions of $n$.

For instance, for $n=3$ we have partitions $(1,1,1)$, $(2,1)$ and $(3)$, and we have conjugacy classes $C_{(1,1,1)} = \{(0)(1)(2)\}$, $C_{(2,1)} = \{(0)(1,2), (0,1)(2), (0,2)(1) \}$, and $C_{(3)} = \{(0,1,2), (0,2,1)\}$. 

Given a partition $\lambda = (\lambda_1, \lambda_2,\ldots, \lambda_k)$ we associate to it the permutation
\[\sigma_{\lambda} = (0, 1, \ldots , \lambda_1-1)(\lambda_1, \lambda_1 + 1, \ldots , \lambda_1+\lambda_2-1)\ldots \]
and call it the {\em canonical representative permutation} of both the partition $\lambda$ and the conjugacy class $C_{\lambda}$ in $S_n$.

Given a set of permutations $\{\sigma_i\} \subset S_n$, the {\em subgroup generated by $\{\sigma_i\}$} is the set of permutations that can be expressed as a composition of elements in $\{\sigma_i\}$ (with repetitions). A subgroup $G$ of $S_n$ is called {\em transitive} if for every pair $i,j \in [n]$ there exists a permutation $\sigma \in G$ such that $\sigma (i) = j$.

Finally, dual to the conjugacy class of a permutation we have its {\em centraliser} $c_{\sigma} = \{ \mu \in S_n \,\mid\, \mu \circ \sigma \circ \mu^{-1} = \sigma   \}$, that is, the set of all elements of $S_n$ that commute with $\sigma$. By construction, we have that $|S_n| = n! = |c_{\sigma}| \cdot |C_{\sigma}|$. Furthermore, if $\{ p_i : a_i \mid 1 \leq i \leq k \}$ is the cycle structure of $\sigma$, then 
\[|c_{\sigma}| = \prod \limits_{i=1}^k a_i! \cdot p_i^{a_i}.\]

\section{The sampling procedure} \label{sampling}

We first give a brief summary of the procedure before giving a detailed discussion and justification of each step.

\begin{algorithm}
    \caption{Sampling graph-encoded surfaces \label{uniform_sampling}}
    \begin{enumerate}
        \item Uniformly sample partition $\lambda$ of $n$, convert to canonical representative $\mu_\lambda$, and uniformly sample permutation $\sigma \in S_n$
        \item Is subgroup generated by $\{ \mu_\lambda, \sigma \}$ transitive? If not, discard $\{ \mu_\lambda, \sigma \}$ and go to step~1
        \item Compute weight $w$ of the graph encoding $(\mu_\lambda, \sigma)$
        \item Add (properties of) $(\mu_\lambda, \sigma)$ to sample with weight~$w$
    \end{enumerate}
\end{algorithm}

\paragraph*{Step 1.} Recall that we always draw a graph-encoded surface with arcs of colour one being vertical, connecting top nodes to bottom nodes. A pair of permutations, say $(\mu,\sigma)$, then describes how arcs of colour two and three connect top to bottom nodes. In this representation, we can still move vertical arcs $(i,i')$ past other vertical arcs without changing the isomorphism type of the graph-encoded surface as a whole. These re-drawings change $\mu$ and $\sigma$, but not their cycle structure. Re-arranging vertical arcs in this way, we can change the graph-encoded surface such its first permutation is the canonical representative of the conjugacy class $\mu$: Just swap $(0,0')$ with the vertical line corresponding to the first letter of the longest cycles of $\mu$, then swap $(1,1')$ with the second letter in this cycle, and so on. A graph-encoded surface written this way is said to be in {\em standard form}. Every graph-encoded surface can be brought into standard form, and a graph-encoded surface in standard form is determined by a partition (giving the second colour), and a permutation (giving the third colour). 

Hence, it is enough to sample a partition $\lambda$ and a permutation $\sigma$ to obtain a sample graph-encoded surface, and every possible graph-encoded surface can be sampled this way with non-zero probability. This drastically reduces the size of the sample space from $n! \cdot n!$ to $\sim e^{\sqrt{n}} \cdot n!$ \cite{partition}. 

To sample a partition uniformly at random, we use the following procedure due to Fristedt \cite{Fristedt1993}: Let $X_1,\dots, X_n$ be $n$ independent geometrically-distributed random variables where for $1\le i\le n$, $X_i$ represents the number of $i$'s in the partition. Also, let $q_i$ be the success-probability parameter for $X_i$, where $q_i = 1 - e^{\frac{-\pi}{\sqrt{6n}}i}$.

We consider $X_i$ to have the support $[0,\infty)$. Thus, the probability mass functions of our random variables are given by $P(X_i = k) = (1 - q_i)^k q_i$.

We sample $X_1,\dots, X_n$ according to their probability mass functions. After sampling, we check if the partition they represent sums up to $n$:

$$\sum_{i=1}^n X_i\cdot i \stackrel{?}{=} n.$$

If this is true, then we have a valid partition. If not, we must re-sample $X_1,\dots, X_n$ and check again.

\begin{example}[Sampling partitions]
    Let $n=7$ and our sample be $X_1 = 3, X_2 = 0, X_3 = 0, X_4 = 1, X_5 = 0, X_6 = 0, X_7 = 0$. It can easily be checked that the sample represents a valid partition (i.e. sums to 7), namely $\lambda = 4+1+1+1$. The canonical representative permutation $\mu_\lambda$ has one cycle of length 4 and three cycles of length 1. To obtain it, simply take the numbers $0\ 1\ 2\ 3\ 4\ 5\ 6$ and insert brackets appropriately to enforce the cycle structure:

    $$\mu_\lambda = (0\ 1\ 2\ 3)(4)(5)(6).$$
\end{example}

Sampling a partition uniformly at random has expected time complexity $O(n^{3/4})$ \cite{Fristedt1993}. Transforming the partition to its canonical permutation can be done in $O(n)$.

To sample a permutation uniformly at random, we shuffle the list $(0,1,\ldots , n-1)$ using Python function \texttt{random.shuffle()} and convert the shuffled list into cycle notation. This has time complexity $O(n)$. This justifies Step 1 of the procedure and shows that the expected running time is $O(n)$.

\paragraph*{Step 2.} We are only interested in connected surfaces. Hence, we need to make sure that the triangulation encoded by $(\mu_\lambda, \sigma)$ is connected. This is a simple test that a graph with $2n$ nodes and $3n$ arcs is connected, and can be done in $O(n)$ time. If a sample is not connected, we discard it and start over. 

\begin{figure}[htb]
    \begin{center}
        \centerline{\includegraphics[width=.48\textwidth]{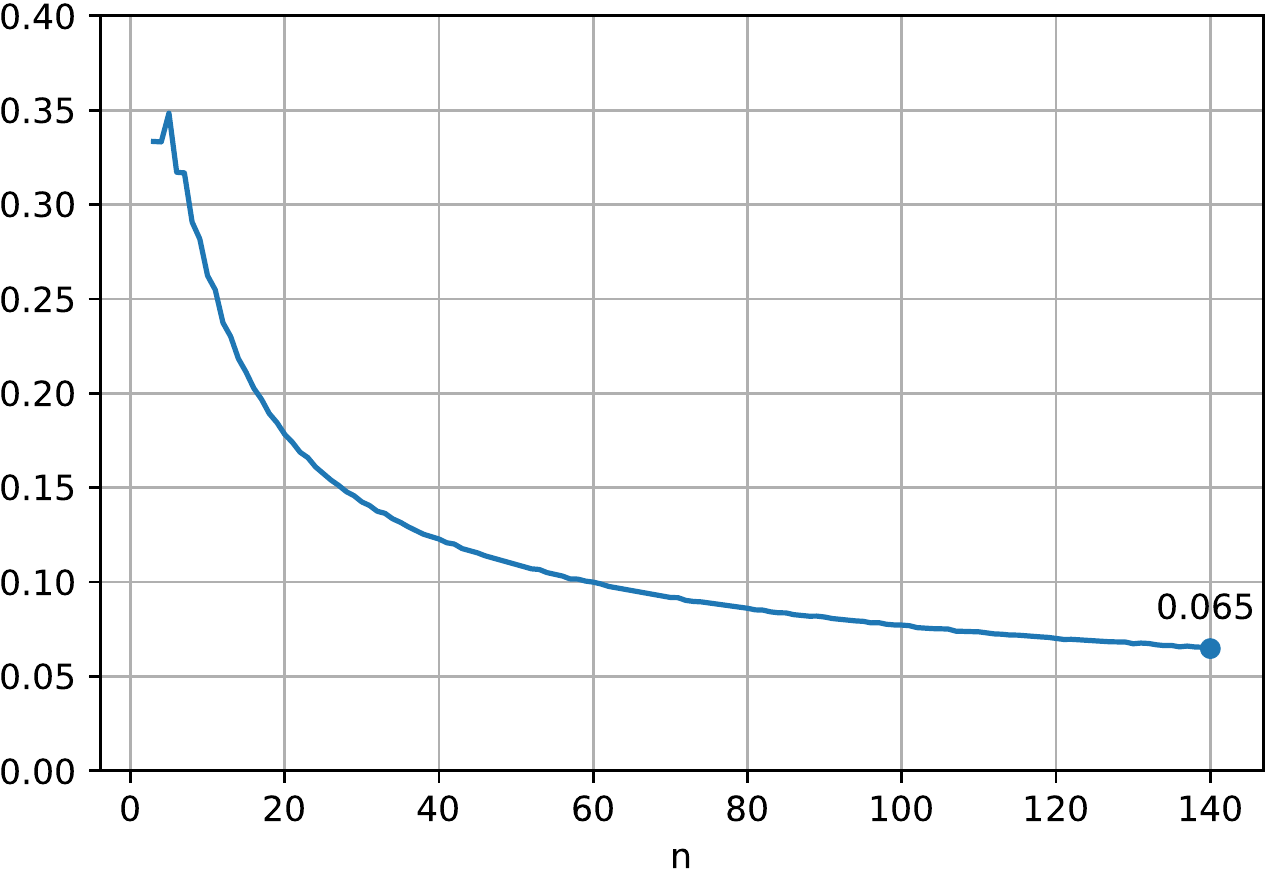}}
    \end{center}
    \vspace{-1.cm}
    \caption{Proportion of graph encodings that are disconnected (and hence discarded) for each $n$.}
    \label{disconnected_gems_proportion}
\end{figure}

The experiments presented in \Cref{experimental_results} show that, for $n \geq 60$, less than $10\%$ of initial samples are discarded due to being disconnected, see \Cref{disconnected_gems_proportion} for more details. Hence, the overall expected running time of Steps 1 and 2 until we obtain a connected graph encoding of a surface can be assumed to be $O(n)$.

\paragraph*{Step 3.} This is the core of the procedure. Given $\mu_\lambda$ and $\sigma$, that is, a graph encoding in standard form, we need to compute the overall number of partition-permutation pairs giving rise to graph encodings in standard form that are isomorphic to $(\mu_\lambda,\sigma)$.

Recall that there are two distinct types of such isomorphisms: those that keep the colours invariant but shuffle around $\sigma$ along symmetries of the bi-coloured cycles defined by $\mu_\lambda$, and those that swap colours but leave everything else invariant. The overall number of symmetries is then the overall number of combinations of these two types of isomorphisms. We handle the two cases separately.

\medskip

First, assume that no colours are swapped. This implies that partition $\lambda$ uniquely determines how the first two colours must be drawn and $\mu_\lambda$ cannot change under the isomorphism. The only way we can redraw $(\mu_\lambda, \sigma)$ to obtain an isomorphic copy $(\mu_\lambda, \tilde{\sigma})$ is by a combination of the following three operations on the bi-coloured cycles of the first and second colour: {\em (a)} Re-orderings of bi-coloured cycles of the same length; {\em (b)} rotations of a single bi-coloured cycle by two arcs; and {\em (c)} a single point reflection of the entire picture through the centre (followed by re-.

If $\mu_\lambda$ has $a_i$ cycles of length $p_i$, $\sum_i a_i \cdot p_i = n$, then there are $a_i!$ ways to re-order those cycles of length $p_i$, and every cycle can be independently rotated in $p_i$ distinct ways. Altogether, this results in 
\[ 2 \cdot \prod \limits_{i=1}^{k} a_i! \cdot p_i^{a_i} = 2 \cdot | c_{\mu_\lambda}|\]
symmetries leaving the first two colours invariant. Here, the extra factor of $2$ comes from the single point reflection around the centre of the drawing.

Note that the cardinality of the centraliser shows up in this formula due to the nature of the transformations we are performing on the graph encoding. This is not a coincidence and the formula can also be deduced purely algebraically.

If $(\mu_\lambda, \sigma)$ has no colour-preserving symmetries, then $2 \cdot | c_{\mu_\lambda}|$ is precisely the number of permutations $\tilde{\sigma}$ such that $(\mu_\lambda, \sigma)$ is isomorphic to $(\mu_\lambda, \tilde{\sigma})$. A colour-preserving symmetry, however, is a symmetry of the first two colours such that $\tilde{\sigma} = \sigma$. Hence, the number of distinct isomorphic copies $(\mu_\lambda, \tilde{\sigma})$ is given by $2 \cdot | c_{\mu_\lambda}| / |\operatorname{Sym} (\mu_\lambda, \sigma)|$, where $\operatorname{Sym}(\mu_\lambda, \sigma)$ denotes the set of colour-preserving symmetries of $(\mu_\lambda,\sigma)$.

\medskip

It remains to look at colour swaps. Since there are three colours, there are $3!=6$ ways to permute them. Given a graph encoding $(\mu_\lambda, \sigma)$ and a permutation of colours, we apply this permutation followed by a re-arrangement of the graph encoding to bring it into standard form. This results in an isomorphic graph encoding $(\mu_{\lambda'}, \sigma')$. If $\lambda' = \lambda$ and $\sigma' = \tilde{\sigma}$ for some $\tilde{\sigma}$ in the list of permutations from the previous step, we have found a colour swap symmetry and this colour swap does not give rise to any additional isomorphic copies. 

If the swapped graph encoding is not equal to one of the previous isomorphic copies, we apply the calculation of colour-preserving isomorphisms from above to $(\mu_{\lambda'}, \sigma')$. We do not have to re-calculate the colour-preserving symmetries, since these remain invariant under isomorphisms (including colour swap isomorphisms).

Note that the set of colour swaps forms a subgroup in the group of colour swaps which is the symmetric group $S_3$. In particular, applying all six colour swaps to $(\mu_\lambda, \sigma)$ detects all colour swap symmetries, and there are either $1$, $2$, $3$ or $6$ such symmetries (the cardinality of a subgroup must divide the cardinality of a group by a classical result from group theory).

\medskip

We can now state our formula for the number of distinct partition-permutation pairs isomorphic to $(\mu_\lambda, \sigma)$ as graph encodings:

$$
\frac{1}{w(\mu_\lambda,\sigma)} = \frac{2}{|\operatorname{Sym}(\mu_\lambda, \sigma)|}
    \sum_{\tau} |c_{\mu_\lambda}|.
$$

This number is the inverse of the weight $w(\mu_\lambda,\sigma)$ applied to $(\mu_\lambda, \sigma)$ before adding it to the sample.

The running time of this step is dominated by the computation of symmetries. A symmetry of a graph encoding is determined by mapping a node and its adjacent three arcs to another node and its adjacent three arcs. There are $2n$ choices for mapping a fixed node, and six choices for mapping its adjacent three arcs. This choice can then be uniquely extended by iteratively mapping neighbours and then their neighbours, etc. The result is either identical to the graph encoding we started with --- in which case we have found a symmetry --- or it is not. Since there are $2n \cdot 6 = 12n$ choices to check, and each choice may require linear running time to be confirmed a symmetry, computing colour-preserving symmetries has a running time of $O(n^2)$. 

Note that, experimentally, an almost linear running time is observed for this step. This is due to the fact that only a negligible number of graph encodings have non-trivial symmetries, see \Cref{avg_symmetries_main} for details. As a result, most initial mappings do not define a symmetry, and in most cases this results in a contradiction to a symmetry after a very small number of steps. 

\paragraph*{Step 4.} Once we have determined a graph encoding together with its weight, we can now perform computations on this graph encoding, or generate a canonical representative of its isomorphism class for inclusion in a weighted uniform sample.

The running time of this step depends on the experiment. If a canonical representative is computed, the procedure of generating this {\em isomorphism signature} is similar in structure to computing the symmetries of a graph encoding and has a worst-case quadratic running time.

\medskip

It follows that the running time to generate one weighted element of our sample splits into two parts: an expected linear running time for all operations up to computing the symmetry, and a worst-case quadratic running time to compute symmetries and isomorphism signature. 

The observed running time in our experiments (where isomorphism signatures were not computed) is overall close to linear; see \Cref{computation_times}.

\begin{example}
  \label{ex:main}
  To better illustrate our sampling procedure, we fully describe it in the case of $n=3$, where numbers are small enough to list the entire sample space.
  
  For $n=3$, we sample a partition $\lambda$ of $3$, and a permutation $\sigma \in S_3$. We thus have $\lambda \in \{ (1,1,1),$ $(2,1),$ $(3) \}$ and $\sigma \in \{ (0)(1)(2),$ $(0,1)(2),$ $(0,2)(1),$ $(0)(1,2),$ $(0,1,2),$ $(0,2,1)\}$ with canonical representative permutations $\mu_{(1,1,1)} = (0)(1)(2)$, $\mu_{(2,1)} = (0,1)(2)$, and $\mu_{(3)} = (0,1,2)$, and hence we get the overall sample space of pairs of permutations presented in \Cref{fig:raw}.
  
  \begin{figure}[htb]
    \renewcommand{\arraystretch}{1.2}
    {\small\begin{tabular}{@{}l@{}||@{}c@{}|@{}c@{}|@{}c@{}|@{}c@{}|@{}c@{}|@{}c@{}}
       &\ (0)(1)(2)\ & \ (0,1)(2)\ &\ (0,2)(1)\ &\ (0)(1,2)\ &\ (0,1,2)\ &\ (0,2,1) \\
       \hline
       \hline
       (0)(1)(2)\ &\raisebox{-0.1cm}{\includegraphics[height=0.4cm]{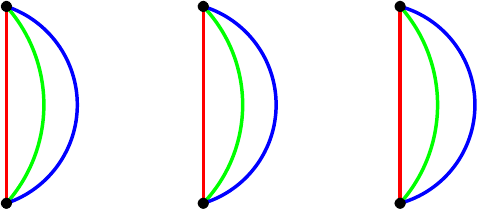}}&\raisebox{-0.1cm}{\includegraphics[height=0.4cm]{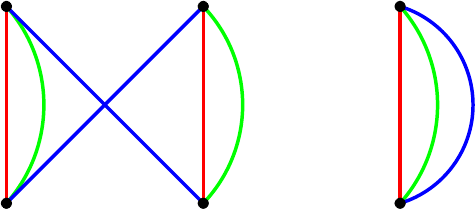}}&\raisebox{-0.1cm}{\includegraphics[height=0.4cm]{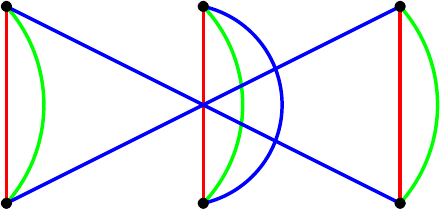}}&\raisebox{-0.1cm}{\includegraphics[height=0.4cm]{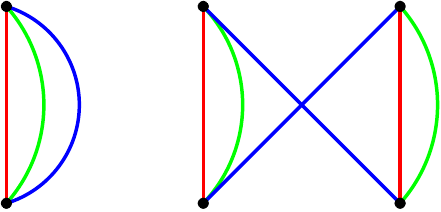}}&\raisebox{-0.1cm}{\includegraphics[height=0.4cm]{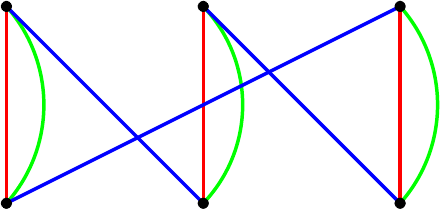}}&\raisebox{-0.1cm}{\includegraphics[height=0.4cm]{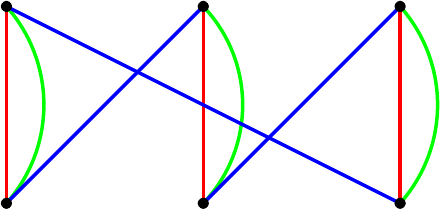}} \\
       \hline
       (0,1)(2)&\raisebox{-0.1cm}{\includegraphics[height=0.4cm]{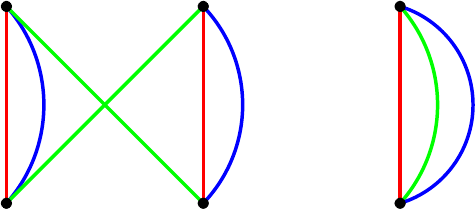}}&\raisebox{-0.1cm}{\includegraphics[height=0.4cm]{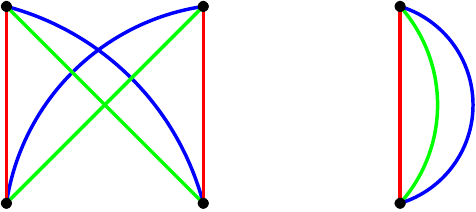}}&\raisebox{-0.1cm}{\includegraphics[height=0.4cm]{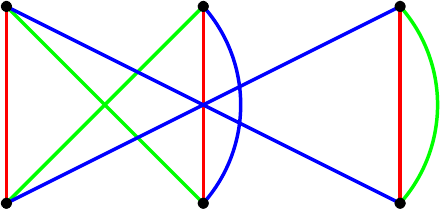}}&\raisebox{-0.1cm}{\includegraphics[height=0.4cm]{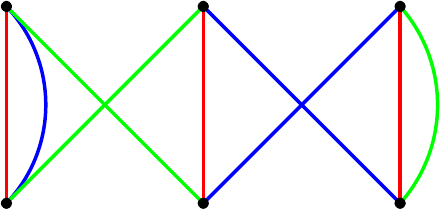}}&\raisebox{-0.1cm}{\includegraphics[height=0.4cm]{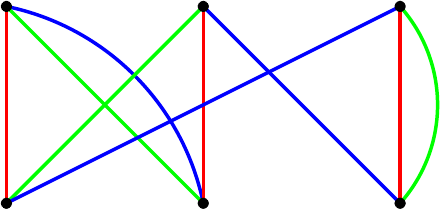}}&\raisebox{-0.1cm}{\includegraphics[height=0.4cm]{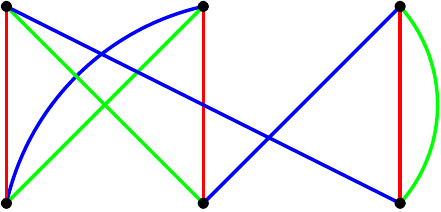}} \\
       \hline
       (0,1,2)&\raisebox{-0.1cm}{\includegraphics[height=0.4cm]{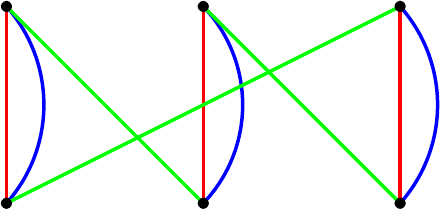}}&\raisebox{-0.1cm}{\includegraphics[height=0.4cm]{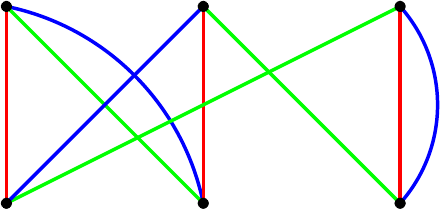}}&\raisebox{-0.1cm}{\includegraphics[height=0.4cm]{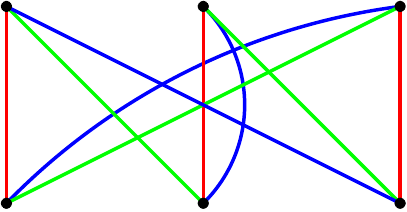}}&\raisebox{-0.1cm}{\includegraphics[height=0.4cm]{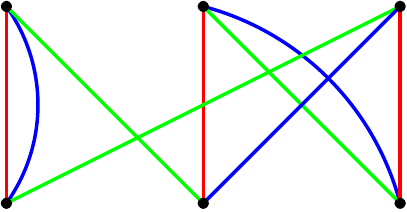}}&\raisebox{-0.1cm}{\includegraphics[height=0.4cm]{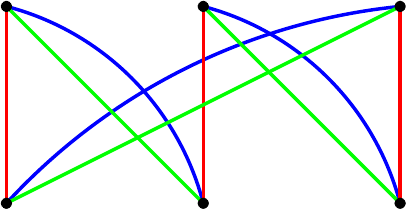}}&\raisebox{-0.1cm}{\includegraphics[height=0.4cm]{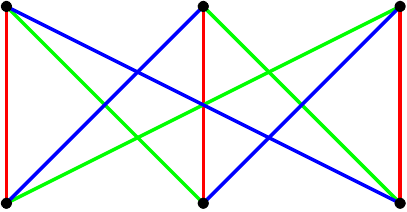}} \\
    \end{tabular}}
    \vspace{-.3cm}  
  \caption{Raw sample space for $n=3$. Disconnected graph encodings are discarded immediately. \label{fig:raw}}
  \end{figure}
  
  After discarding disconnected graph encodings (shaded cells in \Cref{fig:raw}), and running isomorphism and symmetry calculations on the remaining samples, we end up with three isomorphism types of graph encodings: two triangulations of the sphere, denoted $S_0$ and $S_1$, and one triangulation of the torus, denoted by $T_0$, see \Cref{fig:done} for the updated table, and \Cref{fig:list} for a drawing of $S_0$, $S_1$ and $T_0$.
  
    \begin{figure}[htb]
    \renewcommand{\arraystretch}{1.2}
    {\small \begin{tabular}{@{}l@{}||@{}c@{}|@{}c@{}|@{}c@{}|@{}c@{}|@{}c@{}|@{}c@{}}
       &\ (0)(1)(2)\ & \ (0,1)(2)\ &\ (0,2)(1)\ &\ (0)(1,2)\ &\ (0,1,2)\ &\ (0,2,1) \\
       \hline
       \hline
       (0)(1)(2)&$\emptyset$&$\emptyset$&$\emptyset$&$\emptyset$&$S_0$&$S_0$ \\
       \hline
       (0,1)(2)&$\emptyset$&$\emptyset$&$S_1$&$S_1$&$S_1$&$S_1$\\
       \hline
       (0,1,2)&$S_0$&$S_1$&$S_1$&$S_1$&$S_0$&$T_0$\\
    \end{tabular}}
    \vspace{-.3cm}  
   \caption{Isomorphism types of graph encodings for $n=3$. From this table we can directly deduce weights of $1/4$ for $S_0$, $1/7$ for $S_1$ and $1$ for $T_0$. \label{fig:done}}
  \end{figure} 
  
  We perform the calculation of the weight for the example $(\mu_{\lambda},\sigma) = ((0,1)(2),$ $(0,2)(1))$. First note that $(0,1)(2)$ has one cycle of length two and one cycle of length one. We hence have for the number of  symmetries of $\mu_\lambda$: $2\cdot |c_{(0,1)(2)}|=2 \cdot (1! \cdot 2^1)\cdot(1! \cdot 1^1) = 4$. These four symmetries include: the identity, a swapping of the first and second vertical lines, a point reflection, and a combination of the latter two. Swapping the first and second lines produces the isomorphic copy $((0,1)(2),(0)(1,2))$. The same is true for the point reflection, and hence a combination of both of them is a symmetry of the graph encoding. Altogether, we have $2 \cdot  |c_{(0,1)(2)}| / |\operatorname{Sym}((0,1)(2), (0,2)(1))| = 2$.
  
    \begin{figure}
    \centerline{\includegraphics[height=0.75cm]{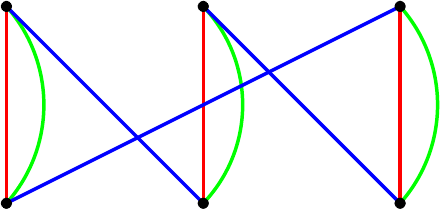} \hspace{0.6cm} \includegraphics[height=0.75cm]{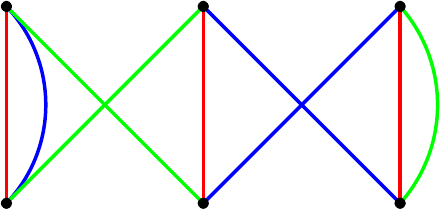} \hspace{0.6cm} \includegraphics[height=0.75cm]{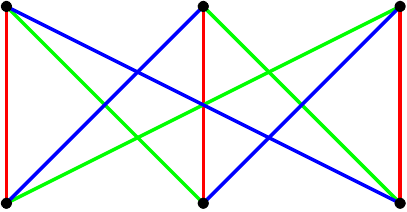}}
  \vspace{-.3cm}
   \caption{Left: $S_0$. Centre: $S_1$. Right: $T$. \label{fig:list}}
  \end{figure} 
  
  Now, swapping colours two and three (followed by rearranging the second colour into standard form) produces the same sample and is thus a colour swap symmetry. Swapping colours one and two produces sample $((0,1)(2), (0,1,2))$ and, again, we have $|c_{(0,1)(2)}|=2$. It follows that swapping colours two and three (a symmetry), followed by swapping colours one and two (which, when combined, is a cyclic shift of colour one to three, three to two, two to one), again results in  $((0,1)(2),(0,1,2))$. The two samples associated with this colour swap are $((0,1)(2),(0,1,2))$ and $((0,1)(2),(0,2,1))$.
  
  Finally, swapping colours one and three or, equivalently, cyclically shifting colours from one to two, two to three, three to one, produces the sample $((0,1,2),(0,2)(1))$. Here, we have $|c_{\mu_{\lambda}}|=3$ and this gives rise to the three samples $((0,1,2),(0,1)(2))$, $((0,1,2),(0,2)(1))$, and $((0,1,2),(0)(1,2))$.
  
  Overall, we have one non-trivial colour swap symmetry and one non-trivial colour preserving symmetry and hence 
  $$
\frac{1}{w} = \frac{2}{|\operatorname{Sym}(\mu_\lambda, \sigma)|}
    \sum_{\tau} |c_{\mu_\lambda}| = \frac{2}{2} (2+2+3) = 7.
$$
  Calculations for $S_0$ and $T$ are analogous. 
\end{example}

\section{Implementation, stability and experiments} \label{experimental_results}

In this section, we present the implementation of our sampling procedure, some calculations on its reliability, as well as the findings of an extended experiment using this implementation.

\medskip

\paragraph*{Implementation.} We use Python 3 as programming language \cite{PythonCoreTeam2020} along with functions from external libraries including \texttt{sympy} and \texttt{numpy} to implement our sampling procedure. We use the \texttt{pandas} and \texttt{matplotlib} libraries to store and plot our data, respectively. The code can be found at the GitHub repository \cite{rajan}.

\medskip

\begin{figure}[htb]
    \begin{center}
        \centerline{\includegraphics[width=.48\textwidth]{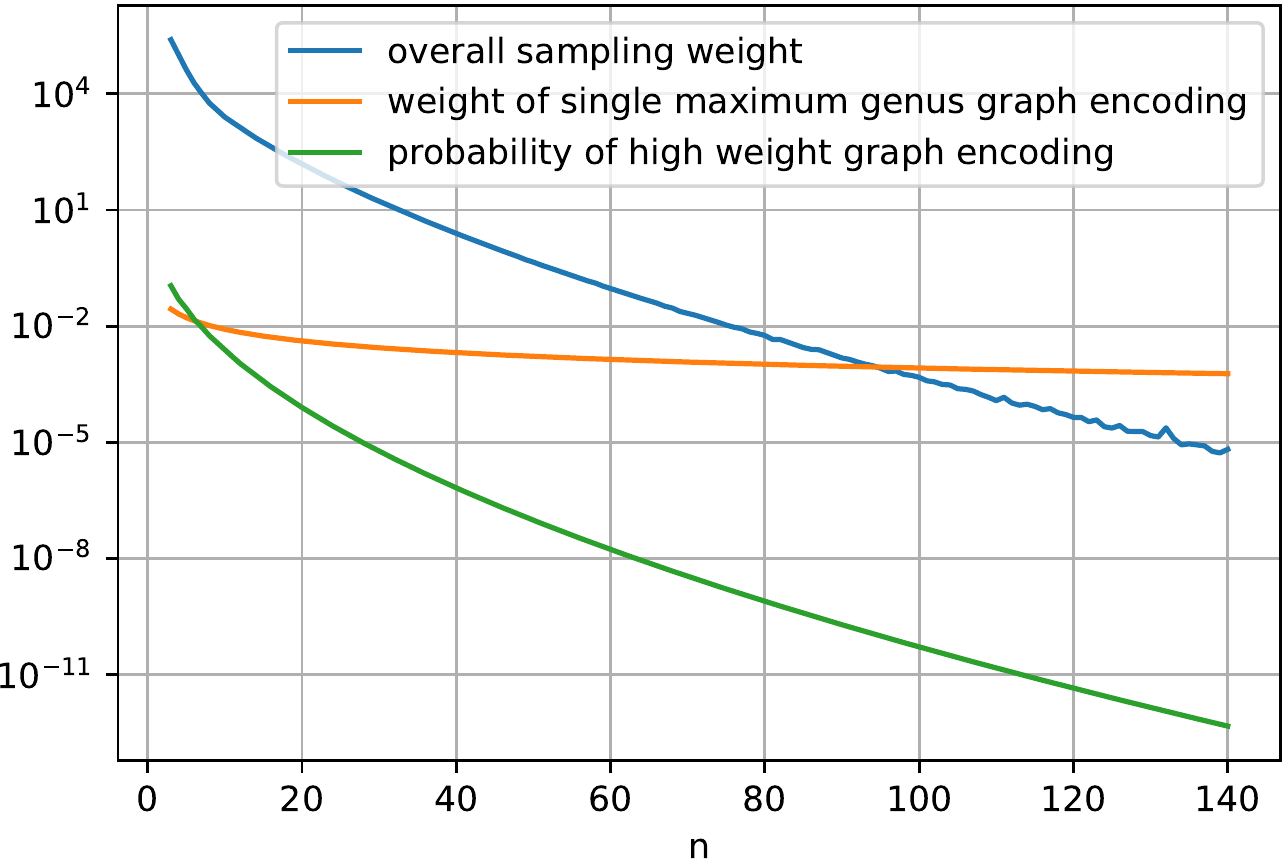}}
    \end{center}
    \vspace{-1.3cm}
    \caption{Comparison of the overall weight of $1,000,000$ samples in our experiment with the weight of a single (non-symmetric) graph encoding of maximum genus.}
    \label{maxweight}
\end{figure}

\paragraph*{Stability.} One issue with our sampling procedure is that weights $w(\mu_{\lambda}, \sigma)$ are very uneven. Specifically, for $\mu_{\lambda}=Id$ the identity we have $|c_{Id}| = n!$ and hence we have $w(\mu_{\lambda}, \sigma) \in O\left(\frac{1}{n!}\right)$. At the same time, weights can be as large as $1$ (see torus $T_0$ in \Cref{ex:main} for an example).  And while this example of a weight of $1$ is very rare, comparatively large weights can occur for large classes of isomorphism classes: As an example, consider a sample where both $\mu_{\lambda} = (1,2, \ldots, n)$ and $\sigma$ only have one cycle. We hence have $|c_{(1,2,\ldots , n)}| = n $ and $w((1,2,\ldots,n),\sigma) \geq \frac{1}{12n}$. These examples coincide exactly with the graph encodings of maximal genus $\lfloor \frac{n-1}{2}\rfloor$ (see the discussion below about Euler characteristic and how to compute the genus of a graph encoding). At the same time, picking such a sample has probability one over the number of partitions of $n$ for $\mu_\lambda$ ($\lambda$ is the unique partition that can give rise to a graph encoding of maximal genus), and $\frac{1}{n}$ for $\sigma$ (there are $(n-1)!$ permutations with one cycle in $S_n$), hence $\sim e^{-\sqrt{n}}/n$ (see \Cref{maxweight} for the exact values for all $n$). 

Since $c_{\mu_{\lambda}}$ becomes very large for permutations with close to $n$ cycles (and hence a low genus), this is not an issue for the relative frequencies of samples with small genus. In \Cref{maxweight} the overall sum of weights of our samples obtained from running our procedure one million times is compared to a single sample of a graph encoding of maximum genus. We conclude from this data that our samples for $n \leq 70$ should still be reliable. At the same time, notice how the quality of our sample visibly reduces for values of $n\geq 100$. We believe that this is a direct consequence of the effect described above.

\begin{figure*}[bt]
    \begin{center}
        \centerline{\includegraphics[width=\textwidth]{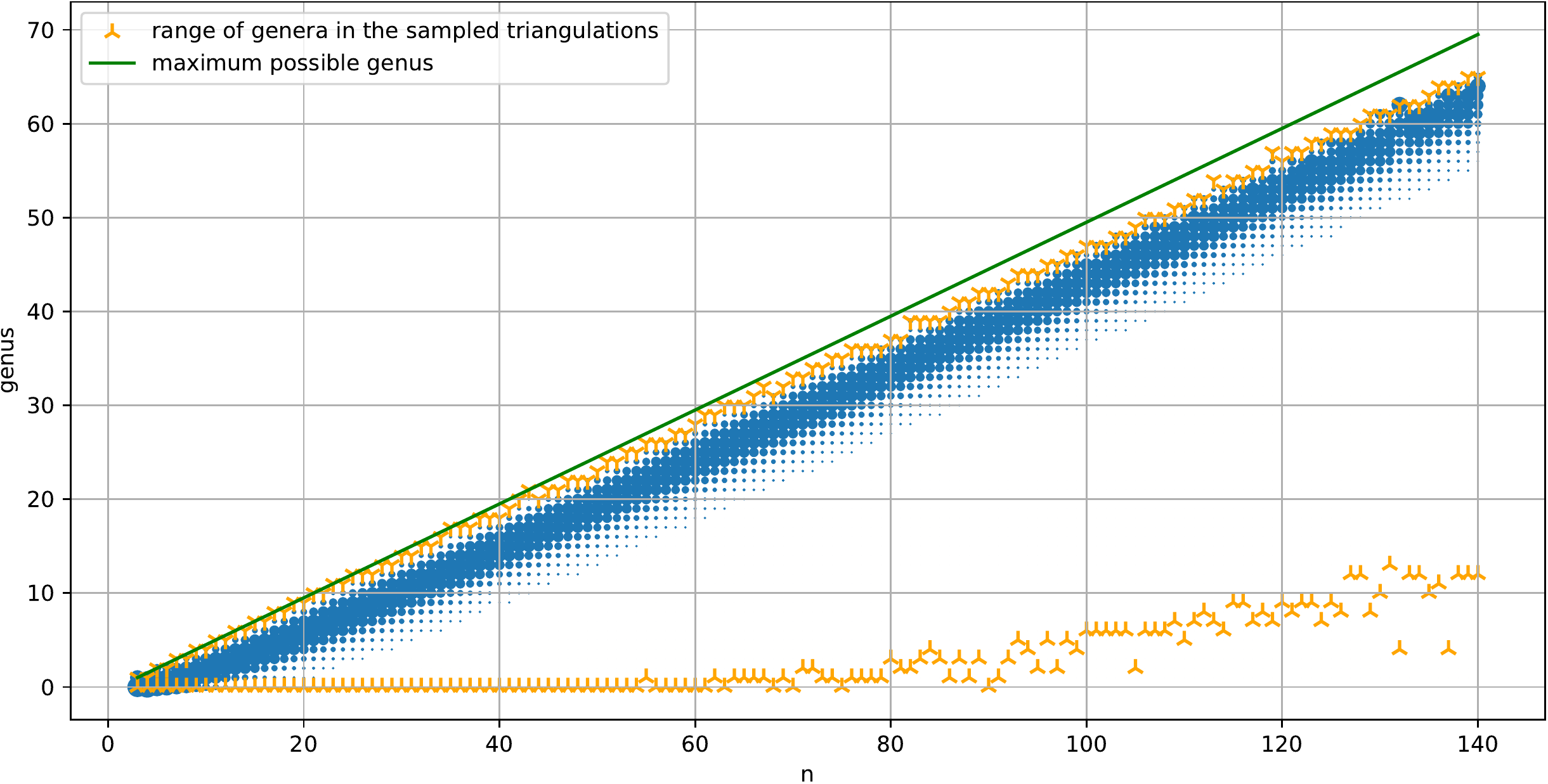}}
    \end{center}
    \vspace{-1.3cm}
    \caption{Empirical distributions of genera for different values of $n$.  Area of circles is proportional to the probability mass of the respective genus in the sample. Notice the distortion for $n=132$ due to a rare sampling of a graph encoding of high genus. This one sample comprises $40\%$ of the entire sample of $1,000,000$ graph encodings.}
    \label{genus_distribution}
\end{figure*}

Finally, note that extremely small weights of size $O\left(\frac{1}{n!}\right)$ must be handled by our procedure. It follows that standard floating point numbers (as used in our procedure) can only be used up to $n \sim 150$. The simple reason for this is that $\frac{1}{n!}$ will be rounded to zero for standard floating point numbers in Python 3 for $n\geq 178$.

\medskip

\paragraph*{The experiment.} We are interested in the following two questions:

\begin{enumerate}
  \item How is the genus distributed in the space of all isomorphism types of graph-encoded surfaces? In particular, what is the expected genus of a uniformly-sampled graph encoding?
  \item For $n$ given, what is the expected number of non-trivial symmetries of a uniformly-sampled graph encoding?
\end{enumerate}

Due to the structure of the sampling procedure, we do not have to do any extra computations to obtain the number of symmetries of a sampled graph encoding.

Similarly, computing the genus from a sampled graph encoding is a straightforward procedure: since the graph encoding is already known to be connected, closed, and orientable, the genus can be computed directly from its {\em Euler characteristic}. The Euler characteristic of a triangulation is given by $\chi = v - e + f$, where $v$, $e$ and $f$ are the number of its vertices, edges and triangles, respectively. We have $f = 2n$ and $e = 3 n$ and hence $\chi = v - n$, and we have for the genus $g = \frac{2 - \chi}{2} = \frac{2 - (v - n)}{2}$. Since we work with permutations rather than triangulations, the number of vertices $v$ is given by the sum of the number of cycles of $\mu_\lambda$, $\sigma$, and $\mu_{\lambda}^{-1} \circ \sigma$. Hence, the running time of this extra step is negligible compared to the running time of the sampling procedure itself.

We ran the procedure one million times for each of the values $3\leq n \leq 140$, and recorded weight, genus, and number of symmetries for each sample.

\medskip

\begin{figure}[htb]
    \begin{center}
        \centerline{\includegraphics[width=.48\textwidth]{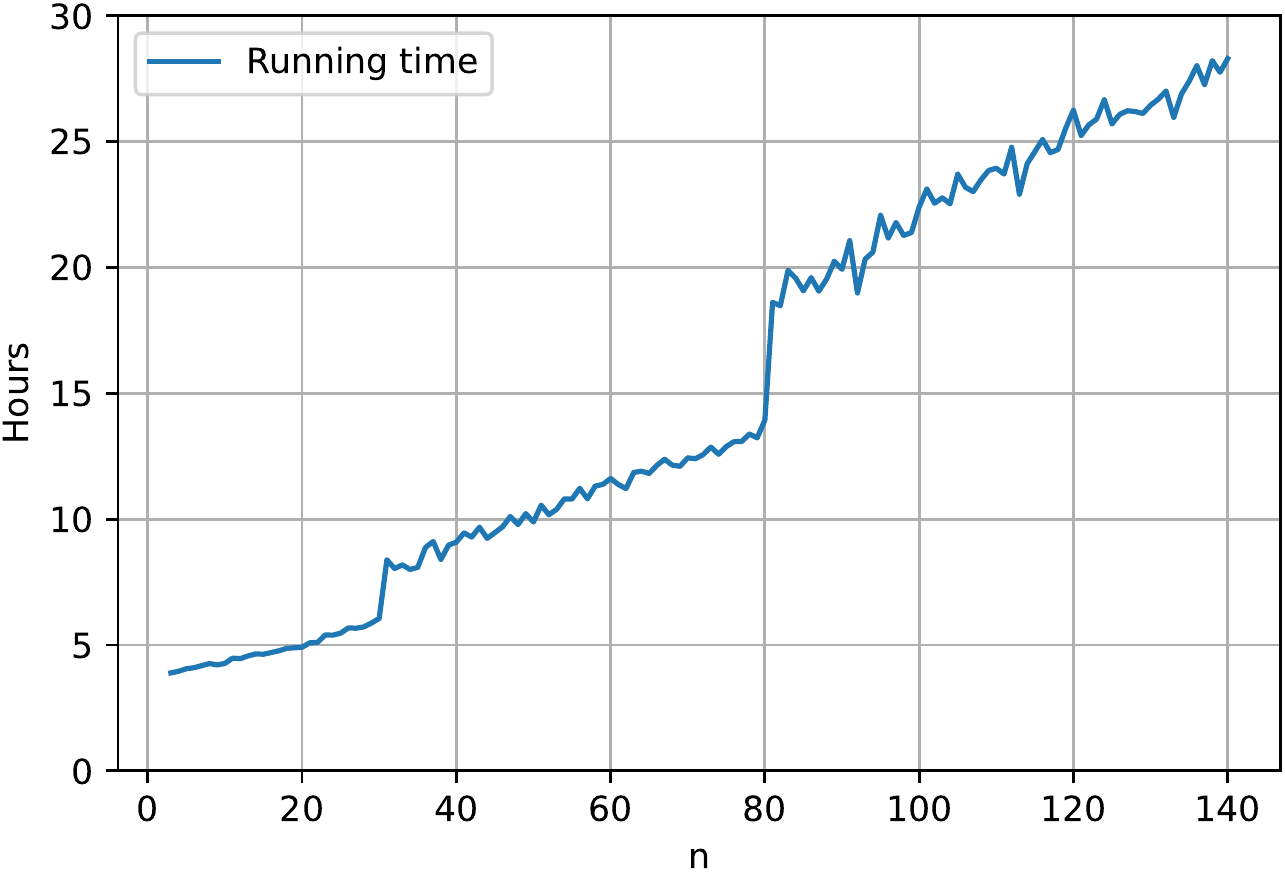}}
    \end{center}
    \vspace{-1.cm}
    \caption{Time required in hours to sample one million triangulations for each value of $n=3,\dots,140$. Note that the jumps are likely artefacts due to the fact that experiments were run in three batches.}
    \label{computation_times}
\end{figure}

\paragraph*{Computation time.} For our computations, we used a computer with $2$ Intel Xeon Gold 6240R CPUs with $2 \times 24$ cores and $192$GB of memory. The procedure is parallelisable without any constraints and we used up to $ 60 $ cores on this machine simultaneously. The observed running time was close to linear in $n$, see \Cref{computation_times} for details.

It follows from this that considerably larger samples as well as samples for larger $n$ can be generated if required.

\medskip

\paragraph*{Genus distribution.} For every value of $n$, we compute the empirical genus distribution of our sample. That is, for every graph encoding in our sample, we record its genus and compute the relative frequencies of all genera in our sample. The result is summarised in \Cref{genus_distribution}. 

Over the entire range of $n$, the bulk of the probability mass of our sample is located on relatively few genera, slightly away from the maximum genus. We want to investigate two questions about the genus distribution: How does the mean genus of a graph encoding with $2n$ triangles behave asymptotically for $n$ tending to infinity? And how concentrated is the probability mass around the mean genus asymptotically? In other words, how does the standard deviation from the mean genus behave, as $n$ tends to infinity?

\medskip

Our experiments, as summarised in \Cref{genus_distribution}, suggest a relatively small gap between the maximum genus and the mean genus. Hence, rather than investigating the mean genus directly, it is more promising to investigate the difference between the maximum genus and the empirical mean genus, raised to different powers to search for linear behaviour.

\begin{figure}[h!]
    \begin{center}
        \centerline{\includegraphics[width=.48\textwidth]{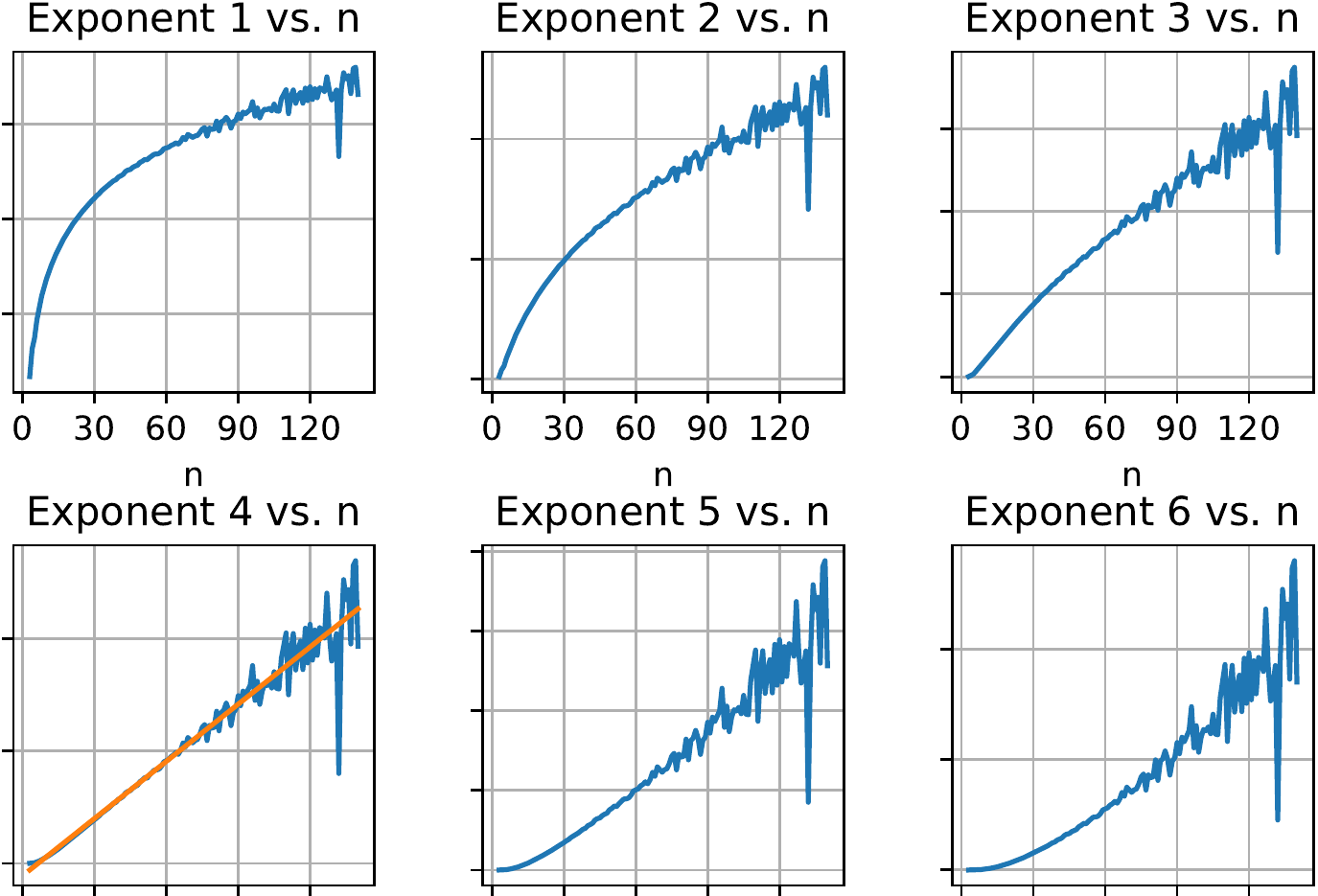}}
    \end{center}
    \vspace{-1.cm}
    \caption{Difference between the maximum genus $\frac{n-1}{2}$ and mean genus raised to powers $1$ to $6$.}
    \label{fig:diff}
\end{figure}

\begin{figure}[hbt]
    \begin{center}
        \centerline{\includegraphics[width=.48\textwidth]{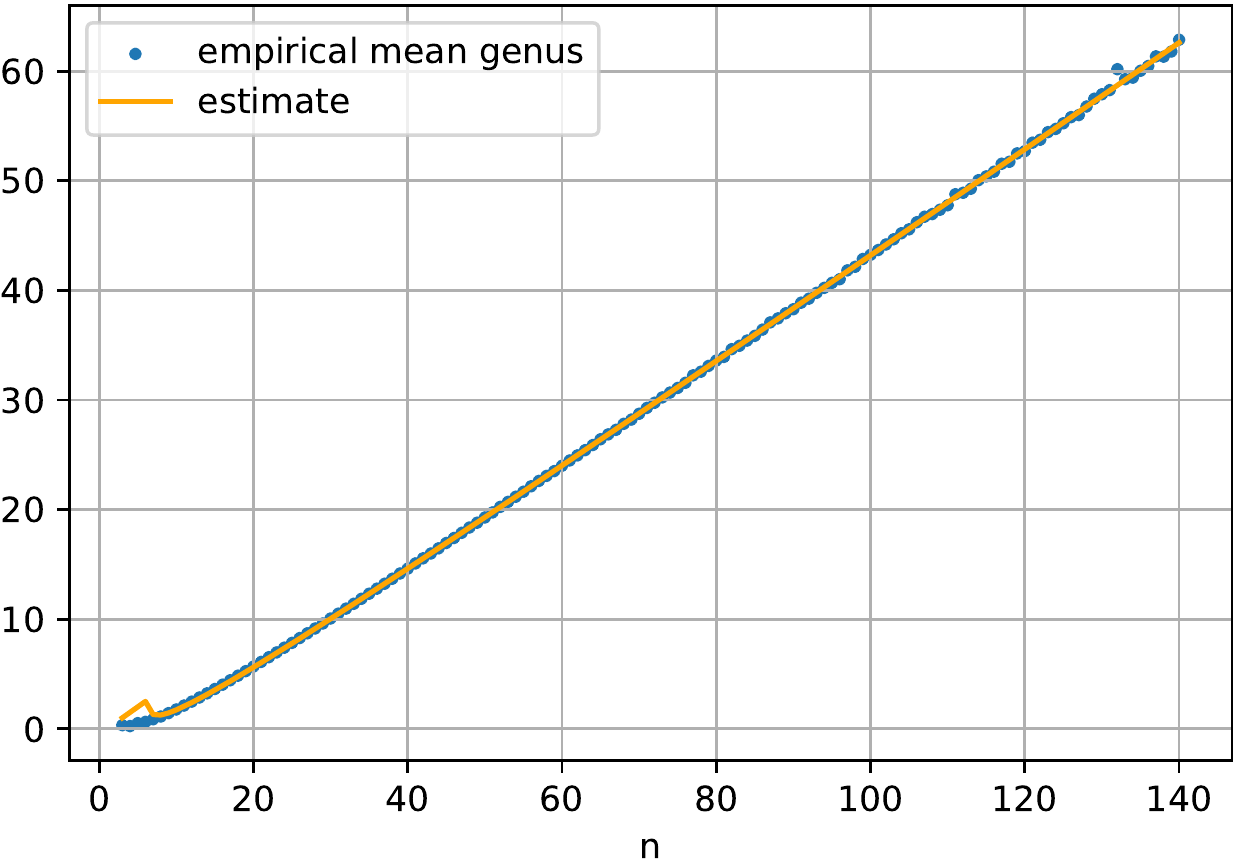}}
    \end{center}
    \vspace{-1.cm}
    \caption{Empirical mean genus compared to our estimate.}
    \label{fig:estimate}
\end{figure}

\Cref{fig:diff} shows the difference between maximum and mean genus raised to the power of the first six positive integers. The difference raised to the power of $4$ exhibits linear behaviour, while lower powers seem sublinear, and higher powers seem superlinear. A regression line $\ell$ through the differences raised to the fourth power has equation $\ell(n) = 16.98n -110.61$. This yields an estimate for the mean genus of $\bar{g} (n) = \frac{n-1}{2} - (16.98n -110.61)^{1/4}$. Note that the estimate was obtained by fitting three parameters (one for the exponent, two for the regression line) to $138$ data points. It is compared with the empirical mean in \Cref{fig:estimate}. We have no theoretical reason to believe that the exponent of $\frac{1}{4}$ is accurate. At the same time, \Cref{fig:diff} suggests that the real asymptotic growth rate of the mean genus will eventually be lower bounded by $\bar{g} (n) = \frac{n-1}{2} - O\left(n^{1/5}\right)$ and upper bounded by $\bar{g} (n) = \frac{n-1}{2} - O \left (n^{1/3}\right )$

\medskip

We can investigate the empirical distribution of the genus further by looking at the changes in standard deviation when varying $n$. 

\begin{figure}[h!]
    \begin{center}
        \centerline{\includegraphics[width=.48\textwidth]{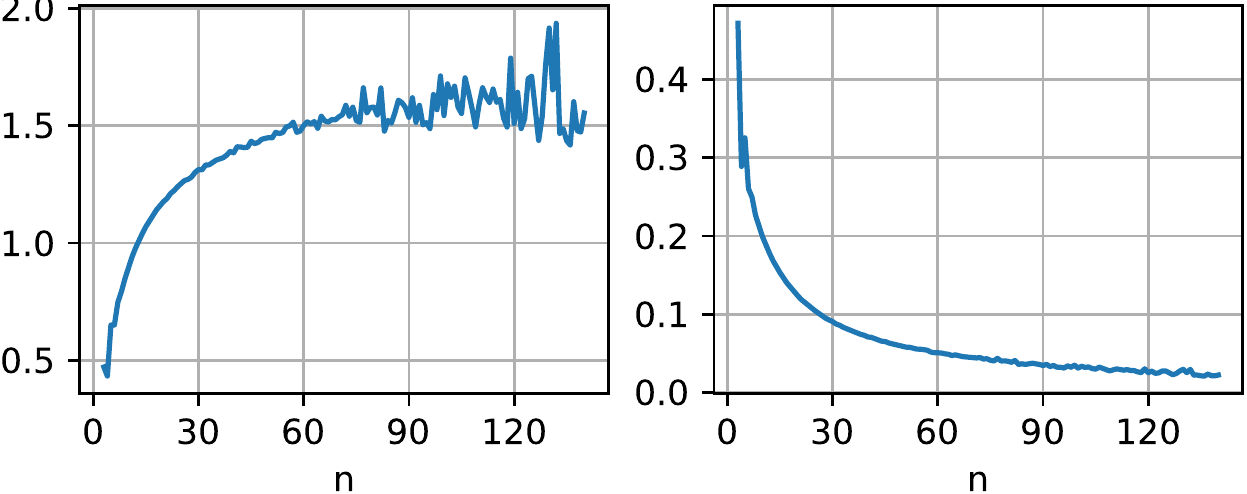}}
    \end{center}
    \vspace{-1.cm}
    \caption{Empirical standard deviation from the mean on scale $0$ to $\lfloor \frac{n-1}{2} \rfloor$ (left) and normalised (right).}
    \label{fig:sd}
\end{figure}

There are essentially two ways of looking at the standard deviation: on a scale from genus $0$ to the maximum genus $\lfloor \frac{n-1}{2} \rfloor$ --- giving us the actual deviation from the mean genus; and on a normalised scale from $0$ to $1$ --- indicating the range of topological types of surfaces we can expect to see compared to the overall number of topological types of surfaces available. Both types of standard deviations are shown in \Cref{fig:sd}. The most obvious finding from this plot is that --- in the normalised setting --- the standard deviation declines with increasing values of $n$. Hence, we see a genus distribution more and more concentrated around the mean genus. Moreover, even in absolute terms the growth of the standard deviation slows down extremely for $n\geq 100$.

\begin{figure}[h!]
    \begin{center}
        \centerline{\includegraphics[width=.48\textwidth]{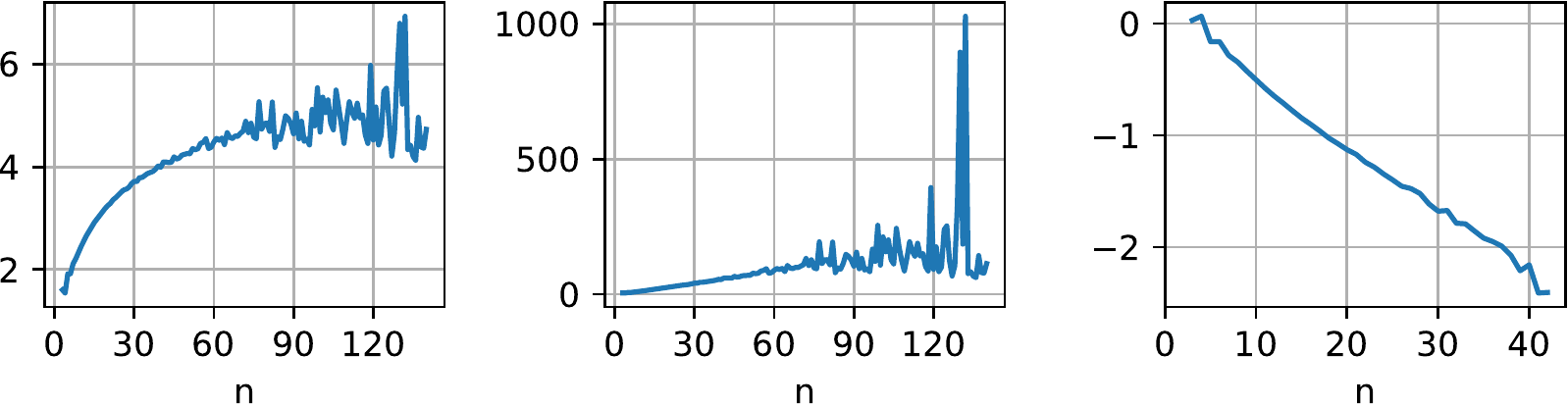}}
    \end{center}
    \vspace{-1.cm}
    \caption{Testing three hypotheses for the growth rate of the absolute standard deviation of the genus from its mean. Left: logarithmic. Centre: doubly logarithmic. Right: exponentially converging to a finite limit. Note that, in the last case, we cannot extend the plot to higher values of $n$ since some data points overshoot the asymptote.}
    \label{fig:hyp}
\end{figure}

To investigate this further, we test three hypotheses for the growth of the standard deviation in absolute terms. The first one is that the standard deviation grows logarithmically, the second one is that it growths doubly logarithmically, and the third one is that it converges to a finite asymptote. For this we apply the exponential and doubly exponential function to the growth rate, and measure the logarithm of the difference of the standard deviation to an asymptote of $1.5$. This is shown in \Cref{fig:hyp}.

Naturally, this experiment suffers much more from the instability (we amplify noise in the standard deviation doubly exponentially and even try to guess an asymptote below certain sample points). Nonetheless, the data seems to be reliable for $n\leq 40$ and provides a clear indicator that the standard deviation grows slower than logarithmically (applying the exponential function is clearly sub-linear, see \Cref{fig:hyp} on the left). Whether or not the standard deviation eventually stops growing or continues to diverge at a very slow rate seems to remain unclear with the data at hand.

Nonetheless, these findings have an important corollary: experimentally this shows that, for $n$ tending to infinity, the genera observed in a uniform sampling will concentrate more and more on a vanishing portion of the overall spectrum of possible genera.

\medskip

\begin{figure}[h!]
    \begin{center}
        \centerline{\includegraphics[width=.23\textwidth]{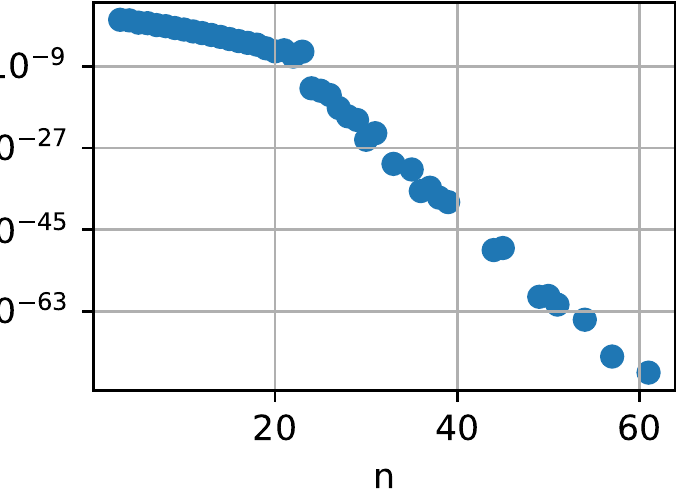}\hfill \includegraphics[width=.23\textwidth]{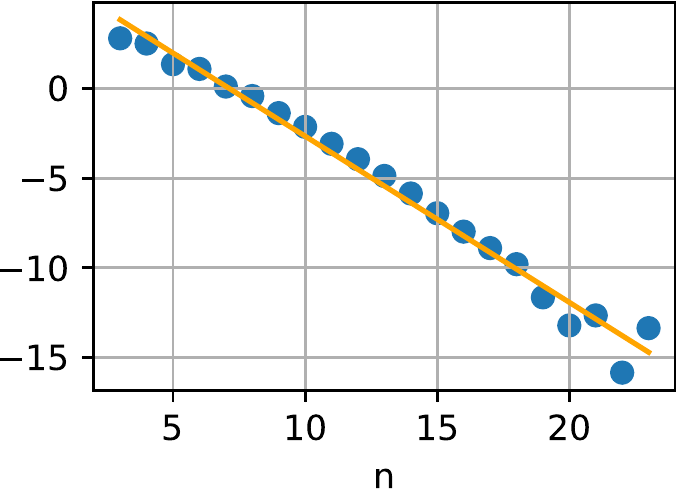}}
    \end{center}
    \vspace{-1.cm}
    \caption{Left: logarithm of the mean number of symmetries minus one for each $n$ for which at least one graph encoding in the sample has a non-trivial symmetry. Right: fitting a line through values of mean symmetries for $n$ between $3$ and $23$ to estimate their rate of decay.}
    \label{avg_symmetries_main}
\end{figure}

\paragraph*{Symmetries.} In this section we investigate the empirical mean for the number of symmetries in our sample. Please note that every graph encoding has one symmetry --- the identity. Moreover, it is expected that for $n$ sufficiently large, most graph encodings do not have any other symmetries. Indeed, \Cref{avg_symmetries_main} on the left reveals that the mean number of non-trivial symmetries rapidly decays as $n$ increases. 

Symmetric graph encodings become incredibly rare as $n$ grows. For $n=32$, not a single graph encoding with non-trivial symmetry is contained in our sample, and the same is the case for all $n\geq 62$. Because of this, it must be assumed that our data is of poor quality, for at least $n\geq 32$. But even under this assumption, the data suggests that the number of non-trivial symmetries decays at a rate super-exponential in $n$: When closing in on the first $20$ data points, for which a reasonable quality can still be assumed, the decay appears to be slightly superlinear on a log scale, see \Cref{avg_symmetries_main} on the right. Fitting a straight line through the first $20$ data points hence suggests that an upper bound for the rate of asymptotic decay is the function  $\exp(-0.95n)$.

\section{Conclusion}

In this report, we presented a uniform sampling procedure for abstract graph encodings of surfaces. We then discussed an implementation of this procedure and an extended experiment spanning $1,000,000$ runs each for numbers of triangles between $6$ to $280$. Within this sample, we then investigated mean genus and mean number of non-trivial symmetries of our samples.

Our main experimental findings are that 

\begin{itemize}
  \item the distribution of genera in a uniform sample seems to concentrate on a vanishing portion of all possible genera for $n$ tending to infinity,
  \item the mean genus seems to be roughly $(16.98n -110.61)^{1/4}$ below the maximum genus possible for any given $n$, and
  \item the mean number of symmetries in a uniform sample likely decays super-exponentially in $n$.
\end{itemize}
  
Finally, we observed that the proportion of graph encodings that were disconnected (and hence discarded) decreased to a relatively small number compared to the size of the overall sample ($< 0.065$).

\medskip

Future work may include validating these findings with larger samples, investigating a theoretical explanation for our observations, proving some of these insights rigorously, and extending the experiments to other properties of graph encodings of surfaces.

Another possible line of future work includes the development of a sampling procedure for graph encodings of a fixed genus.

An additional challenge is to extend this sampling procedure to $3$-dimensional manifolds. Such a procedure would serve as a highly useful tool to investigate $3$-manifolds as well as the performance of algorithms operating on triangulations of $3$-manifolds.

\subsection*{Acknowledgements}

    For this project, the first author received support from the Australian Mathematical Sciences Institute under the Vacation Research Scholarship program. Research of the second author is supported in part under the Australian Research Council’s Discovery funding scheme
(project number DP220102588). The computing power to collect millions of triangulation samples was provided by the School of Mathematics and Statistics at The University of Sydney.

    This project was inspired by the outcomes of several student research projects about graph-encoded surfaces on different levels. The authors would like to acknowledge the members of these projects     
    James Bang,
    Dante Burn,
    Steven Condell,
    Melody Gong,
    Sherry Gong,
    Adam He,
    Chantal Kander,
    Timothy Lapuz,
    Vladimir Levin,
    Harrison Peters,
    Michael Raco,
    Aryonn Rawol,
    Taylor Ruber,
    Zev Shteinman,
    Amelie Skelton,
    Freya Stevens,
    Max Tobin and
    Chloe Yu
    (in strictly alphabetical order) for their hard work and dedication to the topic. The authors would also like to thank Uri Keich for useful discussions about statistics.


\end{document}